\documentclass[12pt]{amsart}
\usepackage{graphics}
\usepackage[all]{xypic}
 \usepackage[mathscr]{eucal}

\newcounter{FNC}[page]
\def\newfootnote#1{{\addtocounter{FNC}{2}$^\fnsymbol{FNC}$%
     \let\thefootnote\relax\footnotetext{$^\fnsymbol{FNC}$#1}}}
\numberwithin{equation}{section}
\theoremstyle{plain}

\newtheorem{thm}[equation]{Theorem}
\newtheorem{lemma}[equation]{Lemma}
\newtheorem{prop}[equation]{Proposition}

\theoremstyle{definition}
\newtheorem{definition}[equation]{Definition}
\theoremstyle{remark}
\newtheorem{remark}[equation]{Remark}

\newcommand{\Der}{\mathrm{Der}}
\newcommand{\Pic}{\mathrm{Pic}}
\newcommand{\spec}{\mathrm{Spec}\;}
\newcommand{\tors}{\mathrm{tors}}
\newcommand{\edim}{\mathrm{edim}}
\newcommand{\bC}{\mathbb{C}}
\newcommand{\bI}{\mathbb{I}}
\newcommand{\bP}{\mathbb{P}}
\newcommand{\bmu}{\boldsymbol{\mu}}
\newcommand{\eC}{\mathscr{C}}
\newcommand{\eK}{\mathscr{K}}
\newcommand{\eM}{\mathscr{M}}

\newcommand{\ch}[2]{\genfrac{(}{)}{0pt}{}{#1}{#2}}
\newcommand{\pe}{\mathbb{P}^2_{E,r}}
\newcommand{\fC}{\mathcal{C}}
\newcommand{\fS}{\mathcal{S}}
\newcommand{\ff}{\mathfrak{f}}

\newcommand{\sC}{\mathcal{C}}
\newcommand{\sE}{\mathcal{E}}
\newcommand{\sI}{\mathcal{I}}
\newcommand{\sL}{\mathcal{L}}
\newcommand{\sN}{\mathcal{N}}
\newcommand{\sO}{\mathcal{O}}
\newcommand{\sP}{\mathcal{P}}
\newcommand{\sU}{\mathcal{U}}

\newcommand{\inertia}{\overline{\sI}_{\bmu}}

\newcommand{\vrho}{\rho}
\renewcommand{\varrho}{{\rho}}

\author{Charles Cadman and Linda Chen}
\title[Enumeration of rational curves tangent to a cubic]{{\bf Enumeration of rational plane curves tangent to a smooth cubic}}
\address{University of Michigan \\
3863 East Hall \\
Ann Arbor, MI 48109-1043}
\email{cdcadman@umich.edu}
\address{Department of Mathematics\\
231 W. 18th Ave.\\
Ohio State University\\
Columbus, OH 43210}
\email{lchen@math.ohio-state.edu}

\begin{document}
\begin{abstract}
We use twisted stable maps to compute the number of rational degree $d$ plane curves having prescribed contacts to a smooth plane cubic.
\end{abstract}
\maketitle

\section{Introduction}

A recent breakthrough in solving enumerative geometry problems in
algebraic geometry occurred in the 1990's with the development of
Gromov-Witten theory and quantum cohomology, inspired by ideas in
physics.  The resulting structures gave surprising and beautiful answers
to classical enumerative geometry problems, for example the enumeration
of degree $d$ rational plane curves passing through $3d-1$ general
points in $\bP^2$.  These structures were defined via moduli spaces of
stable maps to $\bP^2$.

More recently, the Gromov-Witten theory of stacks was developed
and used to define (quantum) orbifold cohomology via moduli
spaces of twisted stable maps to smooth Deligne-Mumford stacks \cite{AGVold,CR}.
In this paper, we apply Gromov-Witten theory
of stacks to enumerative geometry.  Consider the following classical enumerative problem on plane curves.

\begin{quote}
Let $D\subseteq\bP^2$ be a plane curve, and consider tuples 
$\alpha=(\alpha_1,\alpha_2,\ldots)$ and $\beta=(\beta_1,\beta_2,\ldots)$.
How many rational plane curves in $\bP^2$ of degree $d$ meet $D$ at  $\alpha_k$ ``assigned'' points with order of contact $k$ and $\beta_k$ ``unassigned'' points with order of contact $k$ (and pass through $3d-1-\sum (k\alpha_k+(k-1)\beta_k)$
general points), if all contacts with $D$ occur at unibranch points? \end{quote}

The case where all contacts are unassigned and have multiplicity 1 was famously solved by Kontsevich using the Gromov-Witten theory of $\bP^2$.  Caporaso and Harris  solved the case of $D$ a line by using generalized Severi varieties \cite{CH}, and Vakil extended this work to the case of $D$ a conic \cite{Va}.  Using relative Gromov-Witten theory, Gathmann studied the problem when $\alpha=0$: for example, when $d=2, \beta=(0,5)$, i.e. the case of conics five-fold tangent to a plane curve $D$, he obtained a solution given by an explicit polynomial in $\deg D$ \cite{Ga}.
Using Gromov-Witten theory of a particular stack, $\bP^2_{D,2}$, the problem when $D$ is a smooth cubic was solved in the thesis of the first author \cite{thesis}, with $\alpha=(\alpha_1,0,0,\ldots)$ and $\beta = (\beta_1,\beta_2,0,0,\ldots)$.  This was done by using the WDVV equations to obtain several recursions and then showing that the invariants are enumerative.

We give an answer to the problem when $D$ is a smooth cubic for every
$\alpha$ and $\beta$ except for $(\alpha,\beta)=(0,e_{3d})$, thus studying arbitrary higher-order tangency conditions.  Our methods involve the twisted Gromov-Witten theory
of stacks $\bP^2_{D,r}$ for arbitrary integers $r$.  This is of a different
flavor from earlier work of the first author, which studies only the case $r=2$. We denote by  $N_d(\alpha,\beta)$ the number of plane curves of degree
$d$ with contact orders to the smooth cubic $D$ given by the
\emph{contact types} determined by $\alpha$ and $\beta$.

We use several different ideas to compute these numbers.  As in earlier work, we use a WDVV equation to provide relations among the invariants.  This is a little tricky, because every WDVV equation is valid for a fixed value of $r$.  It also turns out that some of the invariants which appear in the equation do not directly count curves, because all of the maps have a component mapping into the cubic with degree $0$.  In order to produce a recursion that involves only enumerative invariants, we relate these nonenumerative invariants to enumerative ones.

It should be noted that the process of showing that these invariants are enumerative is rather involved and cumbersome.  Our approach was inspired by that of Vakil \cite{Va}, in which the irreducible components of the space of stable maps are handled individually.  The point is that, a priori, there can be irreducible components for which the general map is a multiple cover of its image, or even has a reducible source curve.  It must be shown that such components cannot contribute to the invariants, or that if they do, their contribution can be calculated.  We found that such contributions occur only when at most two contacts occur between the rational curve and the cubic.  Moreover, when there are precisely two contacts, the extra contribution can easily be subtracted, so only the case of a single contact remains unsolved by this work (these are the invariants $N_d(0,e_{3d})$).  For small $d$, the contribution to this number coming from a fixed point of order $3d$ was computed by Takahashi \cite{tak_curves}.  He gave a conjectural relation between these numbers and relative invariants in \cite{tak_log}.

In addition to the recursion coming from WDVV, we use the Caporaso-Harris recursion relative to the cubic.  We show that this is valid by proving a relation in the operational Chow ring of the stack of twisted stable maps to $\bP^2_{D,r}$.  This does not give a simple recursion in this case because it only relates invariants for which the degree $d$ and the sequence $\alpha+\beta$ are the same.  However, we prove that all such invariants can be reduced to a single number $M_d(\alpha+\beta)$ via the intriguing formula
$$N_d(\alpha,\beta) = \prod k^{\beta_k}\cdot \left(\sum k\beta_k\right) \cdot \frac{(\sum\beta_k-1)!}{\prod(\beta_k!)}\cdot M_d(\alpha+\beta).$$
These numbers $M_d(\alpha+\beta)$ can be defined independently of the numbers $N_d(\alpha,\beta)$.

While we were not able to compute the numbers $N_d(0,e_{3d})$, we have computed certain nonenumerative analogues, namely the corresponding Gromov-Witten invariants of $\bP^2_{D,r}$ for $r$ sufficiently large relative to $d$.  These numbers seem to agree with the corresponding relative invariants computed by Andreas Gathmann in \cite{Ga2}.  Since both sets of numbers involve virtual contributions, this agreement was unexpected.  It was shown by Maulik and Pandharipande that the relative invariants can be deduced from absolute invariants of $\bP^2$ \cite{maupan}.  We are not aware of any proof of enumerativity for relative invariants.

Finally, we remark that most of the results of this paper hold with the same proof when the plane cubic is replaced with a smooth anti-canonical divisor in a smooth Fano surface.  Then $3d$ should be replaced with the intersection number between the class of the rational curve and the anti-canonical class.  The only change would occur in the WDVV equation (\ref{WDVV_equation}), which was specific to $\bP^2$.

\subsection*{Notation}

Throughout this paper, $E\subseteq\bP^2$ is a fixed smooth plane cubic.  For any positive integer $r$, $\bP^2_{E,r}$ is the $r$th root of $\bP^2$ along $E$, which is constructed in \cite[\S 2.2]{stacks}.  Locally, $\pe$ is the quotient of a cyclic $r$ sheeted covering of an open subset of $\bP^2$, totally ramified along $E$, by the $\mu_r$ action.  We always work over $\bC$ and give irreducible components the reduced induced structure.

If $\alpha = (\alpha_1,\alpha_2,\ldots)$ is a sequence of integers, all but finitely many of which are $0$, we use the following notations (adopted from \cite{CH}).
\begin{equation*}
|\alpha| = \sum_i \alpha_i, \;\;
I\alpha = \sum_i i\alpha_i, \;\;
I^{\alpha} = \prod_i i^{\alpha_i}, \;\;
\alpha! = \prod_i \alpha_i!
\end{equation*}
We use the notation $|\alpha|$ more generally for finite sequences with arbitrary indexing.  Finally, we use $e_i$ to refer to a sequence (finite or infinite) with all zeroes except a $1$ in the $i$th place.

\numberwithin{equation}{subsection}

\section{Morphisms from twisted curves into $\pe$}
\label{sec:def}

In this section we review some facts about twisted stable maps into the $r$th root stack $\pe$ and then study their deformations.  When the domain is smooth and does not map into $E$, a twisted stable map is equivalent to an ordinary map with tangency conditions to $E$ imposed.  Note that the only morphisms from $\bP^1$ to $E$ are constant, a fact which influences the shape of our arguments.  In order to handle the deformation theory of morphisms with contact conditions imposed, we review some results of Caporaso and Harris, which were used in their influential work on enumerating curves with contacts imposed relative to a line.  Then we derive some consequences for the stack of genus $0$ twisted stable maps to $\pe$, which are important to deduce that Gromov-Witten invariants of $\pe$ are enumerative.  Finally, we review some facts about twisted nodes.

\subsection{Review of twisted stable maps}

Twisted stable maps were defined by Abramovich and Vistoli in \cite{AV}.  We are interested in twisted stable maps to $r$th root stacks, so we now recall some results from \cite{stacks}.  First we fix some notation.  There is a smooth divisor $E^{1/r}\subseteq \pe$ which is a gerbe over $E$ banded by $\mu_r$, where $\mu_r$ is the cyclic group of $r$-th roots of unity in $\bC$.  \cite[2.4.4]{stacks}.  Then the following diagram commutes, where $\pi$ exhibits $\bP^2$ as the coarse moduli space of $\pe$.
$$\xymatrix{E^{1/r} \ar@{^{(}->}[r] \ar[d] & \pe \ar[d]^{\pi} \\
E \ar@{^{(}->}[r] & \bP^2}$$
This is not a Cartesian diagram, because $\pi$ is ramified along $E^{1/r}$ to order $r$.  We often use the notation $\sO(\frac{1}{r}E)$ instead of $\sO(E^{1/r})$.  This is consistent with the fact that $\sO(E^{1/r})^{\otimes r}$ is canonically isomorphic to $\pi^*\sO(E)$.  However, we always use the Picard group with integer coefficients, so the fractions should always be interpreted in this way.  This all remains true if $E$ is replaced by any effective Cartier divisor in a scheme, except that $E^{1/r}$ will not be smooth in general.

Recall that a twisted stable map to $\pe$ over a scheme $S$ is a commutative diagram
\begin{equation}
\label{tw_st_map}
\xymatrix{\Sigma_i \ar@{^{(}->}[r] \ar[d] & \fC \ar[r] \ar[d] \ar@/_/[dd] & \pe \ar[d] \\
\sigma_i \ar@{^{(}->}[r] \ar[dr]_{\cong} & C \ar[d] \ar[r] & \bP^2 \\
 & S,}
\end{equation}
where $C\to\bP^2$ is an ordinary stable map over $S$, with sections $\sigma_i$, $\fC$ is a twisted curve with coarse moduli space $C$, $\fC\to\pe$ is representable, and $\Sigma_i$ (the markings) are \'{e}tale gerbes over $\sigma_i$.  These gerbes arise from applying root constructions along $\sigma_i$.  The morphism $\fC\to C$ is an isomorphism away from the gerbes $\Sigma_i$ and the singular locus of $\fC\to S$.  The fibers of $\fC\to S$ can have twisted nodal singularities, which we discuss in section 2.4.

We now recall the notion of \emph{contact type}, which is an integer associated to each marked point $x$ of a twisted stable map, and is a discrete invariant.  When $x$ lies on an irreducible component which does not map into $E$, the contact type equals the multiplicity of the preimage of $E$ at $x$, modulo $r$.  When the component maps into $E$, this is of course no longer true, but a very nice feature of twisted stable maps is that such an invariant still exists, and comes quite naturally out of the theory.  See \cite[\S 3.3]{stacks} for details.  For an $n$-pointed stable map, the contact type forms an $n$-tuple $\vrho=(\varrho_1,\ldots, \varrho_n)$, where $0\le \varrho_i\le r-1$.

The following proposition follows from \cite[3.3.6]{stacks} and the fact that any rational curve mapping into $E$ must do so with degree $0$.

\begin{prop}
\label{contact_type}
Let $\fC$ be a smooth, $n$-pointed, genus $0$ twisted curve over a scheme $S$ and let $\ff:\fC\to\pe$ be a twisted stable map of positive degree and contact type $\vrho$.  Let $C$ be the coarse moduli space of $\fC$ with induced markings $\sigma_i\subseteq C$, and let $f:C\to\bP^2$ be induced by $\ff$.  Then there is an effective Cartier divisor $Z\subseteq C$ such that
\begin{equation}
\label{eq:contact_cond}
f^*D=rZ+\sum_{i=1}^n\varrho_i\sigma_i.
\end{equation}
Moreover, given a morphism $f:C\to\bP^2$ and an effective Cartier divisor $Z\subseteq C$, there is a unique (up to isomorphism) twisted curve $\fC$ with coarse moduli space $C$ and a unique twisted stable map $\ff:\fC\to\pe$ with contact type $\vrho$ which induces $f$.
\end{prop}

The contact type is defined as follows.  Recall that $\Sigma_i$ is obtained from $\sigma_i$ by a root construction.  Let $r_i'$ be the order of this root.  Given a morphism $\ff:\fC\to\pe$, where $\fC$ is a twisted curve over a point, the restriction of the line bundle $\ff^*\sO(\frac{1}{r}E)$ to $\Sigma_i$ determines a character of $\mu_{r_i'}$, since $\Sigma_i\cong B\mu_{r_i'}$.  This isomorphism is fixed so that the tangent space at $\Sigma_i$ corresponds to the standard inclusion $\mu_{r_i'}\subseteq\bC^*$.  From the character induced by $\ff^*\sO(\frac{1}{r}E)$, we obtain a unique integer $k_i$ between $0$ and $r_i'-1$ such that this character is the $k_i$-th power of the standard inclusion.  Representability of $\ff$ implies that $\gcd(k_i,r_i')=1$, and in fact this is equivalent to $\ff$ being representable along $\Sigma_i$.  The contact type $\vrho$ is defined by $\varrho_i = k_ir/r_i'$.  
Note that $r_i'$ and $k_i$ can be recovered from $\varrho_i$ by the formulas
\begin{equation}
\label{rk_from_rho}
r_i' = \frac{r}{\gcd(r,\varrho_i)},\;\;\; k_i = \frac{\varrho_i}{\gcd(r,\varrho_i)}.
\end{equation}

\begin{lemma}
\label{contact_deg0}
Fix an embedding $B\mu_r\subseteq\pe$.  If $\ff:\fC\to B\mu_r\subseteq\pe$ is a representable morphism, where $\fC$ is a twisted curve with coarse moduli space $\bP^1$, then the sum of the contact types at the marked points of $\fC$ is divisible by $r$.  Conversely, given $n$ distinct points $x_1,\ldots,x_n\in\bP^1$ and an $n$-tuple $\vrho$ of integers with $0\le\varrho_i\le r-1$, such that $\sum\varrho_i$ is a multiple of $r$, there is a unique (up to isomorphism) twisted curve $\fC$ with coarse moduli space $\bP^1$ and representable morphism $\ff:\fC\to B\mu_r$ with contact types given by $\vrho$.
\end{lemma}

\begin{proof}
Given a morphism $\ff:\fC\to B\mu_r$ with these contact types,  $$\ff^*\sO(E^{1/r})\cong \sO\left(\sum_{i=1}^n \frac{\varrho_i}{r}x_i - a\right)$$ for some integer $a$.  Since the $r$-th power of this line bundle is trivial, it follows that $\sum_{i=1}^n \varrho_i = ar.$  Since this is the only degree $0$ line bundle on $\fC$ with these contact types, the morphism $\ff$ is determined by $\vrho$.  Moreover, $\fC$ is determined by $\vrho$ in light of (\ref{rk_from_rho}).
\end{proof}

Now we introduce some more notation.  Since we are only interested in genus $0$ stable maps to $\pe$, we use $\eK_d^r(\vrho)$ to refer to the stack $\eK_{0,n}(\pe,d,\vrho)$, where $n$ is the dimension of the vector $\vrho$, and $$\eK_{0,n}(\pe,d,\vrho)\subseteq\eK_{0,n}(\pe,d)$$ is the open and closed substack of the stack of $n$-pointed, genus $0$, degree $d$ stable maps to $\pe$ having contact type $\vrho$.

We often use the following variant of this notation.  Let $\gamma = (\gamma_0,\ldots,\gamma_{r-1})$ be an $r$-tuple of integers.  From this, we construct a $|\gamma|$-tuple $\vrho$ by the rule
$$\varrho_i=j\mathrm{\ if\ } \sum_{k=0}^{j-1}\gamma_k < i \le \sum_{k=0}^j \gamma_k.$$
In other words, we want precisely $\gamma_j$ of the integers $\varrho_i$ to equal $j$.  We sometimes write $\eK_d^{\gamma}$ for $\eK_d^r(\vrho)$, and alternate between these two notations.  We sometimes abuse notation and write $\eK_d^{\gamma}$ when $\gamma=(\gamma_1,\gamma_2,\ldots)$ is an infinite sequence and $r$ is an arbitrary integer such that $\gamma_i = 0$ for $i\ge r$.  In this case, the stacks $\eK_d^{\gamma}$ depend on the suppressed integer $r$, but this dependence is not very significant for $r$ large enough.  We implicitly set $\gamma_0=0$ in this case, which means that all markings are twisted.

The expected dimension of the stack $\eK_d^r(\vrho)$ is
\begin{equation}
\label{eq:edim}
\frac{1}{r}(3d - \sum \varrho_i) + n - 1,
\end{equation}
where $n$ is the number of entries in $\vrho$.  If $\vrho$ comes from an $r$-tuple $\gamma$ as above, then this equals $$\frac{1}{r}(3d-I\gamma) + |\gamma| - 1.$$  This formula is a special case of equation 3.5.1 of \cite{GWinvs}.  

Recall that the expected dimension is the degree (under the grading by dimension) of the virtual fundamental class in the Chow group of $\eK_d^r(\vrho)$.  The virtual fundamental class appears in Gromov-Witten theory as the class which one integrates against, as in Definition \ref{GWI_def}.  The expected dimension is always less than or equal to the actual dimension of any irreducible component, with equality if and only if the virtual fundamental class equals the usual fundamental class on that component.  These facts are true whenever one has a perfect obstruction theory \cite{BF}.  In particular, this is true in twisted Gromov-Witten theory \cite[\S 4.5]{AGV}, \cite[\S 3.1]{GWinvs}.

%%\begin{definition}
%%\label{nice_maps}
%%We define $\eK^{*}_{0,n}(\pe,d,\vrho)$ as follows.  Define an irreducible component of $\eK_{0,n}(\pe,d,\vrho)$ to be \emph{good} if the general map in that component has a smooth source curve and maps birationally onto its image.  Let $\eU\subset\eK_{0,n}(\pe,d,\vrho)$ be the open substack obtained by removing all the bad irreducible components, and let $\eK^*_{0,n}(\pe,d,\vrho)$ be its stack-theoretic closure.
%%\end{definition}

%%Note that $\eK^*_{0,n}(\pe,d,\vrho)$ has an open dense substack which is a scheme.  Indeed, the general map in each irreducible component has no automorphisms, and $\eK_{0,n}(\pe,d)$ has a projective coarse moduli scheme \cite[1.4.1]{AV}.

\subsection{Maps from smooth twisted curves}

We begin with a lemma about infinitesimal deformations preserving contacts.  Let $C$ be a smooth curve (as always, over $\bC$) and let $X$ be a smooth variety.  Let $f:C\to X$ be a morphism, let $D\subseteq X$ be a smooth divisor, and let $c_1,\ldots, c_n\in C$ be distinct points such that
\begin{equation}
\label{divisor_equality}
f^*D=\sum_{i=1}^n m_ic_i
\end{equation}
for some positive integers $m_i$.  We assume that (\ref{divisor_equality}) is an equality of \emph{subschemes} of $C$; in particular, $f$ is nonconstant.  Let $\sN$ be the cokernel of the differential $$T_C\to f^*T_X.$$  Let $Def^1_{log}(f)$ denote the space of first order infinitesimal deformations of the tuple $$(f:C\to X,c_1,\ldots, c_n),$$ with $X$ fixed, which preserve condition (\ref{divisor_equality}).

\begin{lemma}
\label{inf_def}
There is a surjective morphism of sheaves
\begin{equation}
\label{morph_inf_def}
\sN\to \bigoplus_{i=1}^n \sO_{(m_i-1)c_i}
\end{equation}
and a canonical isomorphism
$$Def^1_{log}(f)\cong H^0(C,\sN^{log})$$\
where $\sN^{log}$ is the kernel of (\ref{morph_inf_def}).  Here $\sO_{mc}$ refers to the structure sheaf of a multiple point of multiplicity $m$ at $c\in C$.
\end{lemma}

\noindent{\it Remark.}  Such deformations are the same as deformations of log morphisms, so this follows from log deformation theory.  We include a proof for completeness.

\begin{proof}  First we recall the construction which identifies first order deformations of $f$ fixing $X$ with $H^0(C,\sN)$.  Let $\bI=\spec \bC[\epsilon]/(\epsilon^2)$, and suppose we have the following commutative diagram.
$$\xymatrix{
 & X \\
C \ar[ur]^f \ar@{^{(}->}[r] \ar[d] \ar@{}[dr]|{\Box} & \sC \ar[d] \ar[u]_F \\
\spec \bC \ar@{^{(}->}[r] & \bI}$$
Cover $C$ by affines $U_i$.  Since $C$ and $\sC$ have the same underlying topological space, we get an open covering of $\sC$ by subschemes $\sU_i$ so that $U_i$ embeds into $\sU_i$ as $(\sU_i)_{\mathrm{red}}$.  Since nonsingular affine varieties have no nontrivial first order deformations, there are isomorphisms $\varphi_i:U_i\times\bI\to\sU_i$.  On the overlaps $U_{ij}$, $\varphi_j^{-1}\circ\varphi_i$ determines a derivation $\mu_{ij}\in H^0(U_{ij},T_{U_{ij}})$.  These form a 1-cycle, and hence determine an element $\mu\in H^1(C,T_{C})$.

On $U_i$, the morphisms $F$ and $\varphi_i$ determine a morphism $\psi_i:U_i\times\bI\to X$, and hence a derivation $\nu_i\in H^0(C,f^*T_X\vert_{U_i})$.  The commutative diagram
$$\xymatrix{
U_{ij}\times\bI \ar[r]^{\varphi_i} \ar[dr]_{\psi_i} & \sU_{ij} \ar[r]^{\varphi_j^{-1}} & U_{ij}\times\bI \ar[dl]^{\psi_j} \\
 & X &}$$
shows that on $U_{ij}$ we have $\nu_i=df(\mu_{ij})+\nu_j$.  It follows that the $\nu_i$ glue to give an element $\nu\in H^0(C,\sN)$.  This is independent of the choices and identifies the first order deformations with $H^0(C,\sN)$ \cite[p. 96]{HM}.

We need a necessary and sufficient condition for the first order deformation corresponding to $\nu$ to preserve the multiplicities of $f^*D$.  For this we can reduce to the affine situation $U_i=\spec S\to \spec R\subseteq X$ and assume that $f^*D=mp$ for some $p\in U_i$.  Let $r\in R$ be a local equation for $D$ and $s\in S$ a local equation for $p$ so that $f^*r=us^m$ with $u\in S$ a unit.  Under the morphism $U_i\times\bI\to X$ determined by $\nu_i\in\Der_{\bC}(R,S)$, $r$ pulls back to $f^*r+\nu_i(r)\epsilon$.  This defines a multiplicity $m$ divisor if and only if there are elements $v,w,x,y\in S$ with $v$ a unit such that $$us^m+\nu_i(r)\epsilon = (v+w\epsilon)(x+y\epsilon)^m.$$  This is equivalent to $\nu_i(r)$ being an element of the ideal generated by $s^{m-1}$.  We now interpret this condition.

By standard arguments, there are maps making the following diagram commutative with exact rows and columns.  Here $\sE$ is the sheaf of logarithmic vector fields on $X$ relative to $D$.
$$\xymatrix{
& 0\ar[d] & 0\ar[d] & 0\ar[d] & \\
0\ar[r] & T_C(-\sum c_i) \ar[d]\ar[r] & f^*\sE \ar[r]\ar[d] & \sN^{log} \ar[d]\ar[r] & 0 \\
0\ar[r] & T_C \ar[d]\ar[r] & f^*T_X \ar[d]\ar[r] & \sN \ar[d]\ar[r] & 0 \\
0\ar[r] & \bigoplus_i \sO_{c_i} \ar[d]\ar[r] & f^*N_{D/X} \ar[d]\ar[r] & \bigoplus_i \sO_{(m_i-1)c_i} \ar[r]\ar[d] & 0 \\
& 0 & 0 & 0 & }$$
The morphism $f^*T_X\to f^*N_{D/X}$ locally sends a derivation $\nu_i$ to $\nu_i(r)$ modulo $s^{m_i}$.  Therefore, the condition that $\nu_i(r)$ is in the ideal generated by $s^{m_i-1}$ means precisely that $\nu_i$ is a section of $\sN^{log}|_{U_i}$.
\end{proof}

\begin{remark}
\label{no_torsion}
The sheaf $\sN^{log}$ constructed above has no torsion supported in $f^{-1}(D)$.  This can be verified locally, and we give the argument when $X$ is a surface.  Suppose that $x,y$ are coordinates on $X$, with $D$ defined by $y=0$.  Let $f:C\to X$ be defined by $$t\mapsto (p(t),q(t)),$$ with $p(0)=q(0)=0$.  Then under the differential, $$t\frac{\partial}{\partial t} \mapsto t p'(t)\frac{\partial}{\partial x} + \frac{tq'(t)}{q(t)}y\frac{\partial}{\partial y}.$$  Since $tq'(t)/q(t)$ is nonvanishing at $t=0$, it follows that the fiber of $\sN^{log}$ at the point $t=0$ is generated by $\partial/\partial x$.
\end{remark}

Proposition \ref{contact_type} implies that the deformation theory of twisted stable maps $\fC\to\pe$ from smooth twisted genus $0$ curves $\fC$ is equivalent to the deformation theory of maps $C\to\bP^2$ from smooth rational marked curves $C$ with contact conditions imposed at the markings.  This explains the importance of the above lemma, as well as the following two results of Caporaso and Harris (stated here for the special case of $E\subseteq\bP^2$).

Let $\pi:C\to B$ be a smooth, proper family of connected curves over a smooth base $B$, let $f:C\to\bP^2$ be a morphism, and let $b\in B$ be a general point.  Assume that no fiber of $\pi$ maps to a point under $f$.  Let $\sN_b$ be the cokernel of the differential $df_b:T_{C_b}\to f_b^*T_{\bP^2}$, which is injective by hypothesis.  We have a morphism $\kappa_b:T_b B\to H^0(C_b,\sN_b)$ induced by the family of morphisms $C\to B\times\bP^2$, and this is often called the Horikawa map in recognition of Horikawa's foundational work \cite{Ho}.

\begin{lemma}
\label{Cap_Harris_1}
{\bfseries \cite[2.3]{CH}}  Let $b\in B$ be a general point and assume that $f_b$ maps $C_b$ birationally onto its image.  Then $$\mathrm{Im}(\kappa_b)\cap H^0(C_b,(\sN_b)_\tors)=0.$$
\end{lemma}

Let $Q\subseteq C$ be the image of a section of $\pi$ such that $f^*E$ has multiplicity $m$ along $Q$.  Let $q=Q\cap C_b$.  Let $\ell-1$ be the order of vanishing of $df_b$ at $q$.

\begin{lemma}
\label{Cap_Harris_2}
{\bfseries \cite[2.6]{CH}}  For any $v\in T_bB$, the image of $\kappa_b(v)$ in $H^0(C_b,\sN_b/(\sN_b)_\tors)$ vanishes to at least order $m-\ell$ at $q$ and cannot vanish to any order $k$ with $m-\ell<k<m$.  If $f(Q)$ is a point, then it vanishes to order at least $m$ at $q$.
\end{lemma}

Note that in the genus $0$ case, there cannot be gaps in the orders of vanishing of sections of a line bundle.  So the previous lemma implies that if $\ell>1$, then the image of $\kappa_b(v)$ either vanishes precisely to order $m-\ell$ for every $v\in T_bB$ or else vanishes at least to order $m$ for every $v\in T_bB$.

\begin{thm}
\label{gen_smooth}
Let $\gamma=(\gamma_1,\gamma_2,\ldots)$ be a sequence of integers such that $I\gamma = 3d$ for some positive integer $d$, and choose $r>3d$.  Let $V\subseteq \eK_d^{\gamma}$ be an irreducible component such that the general map parametrized by $V$ is from a smooth curve mapping birationally onto its image.  If $|\gamma|\ge 2$, then 
\begin{enumerate}
	\item the dimension of $V$ equals the expected dimension, $|\gamma|-1$,
	\item $\eK_d^{\gamma}$ is generically reduced along $V$, and
	\item every evaluation map $\eK_d^{\gamma}\to E$ corresponding to a twisted marking is surjective.
\end{enumerate}
Moreover, if $f:C\to\bP^2$ is the map of coarse moduli spaces corresponding to a general point of $V$, then the following are true.
\begin{enumerate}
	\item[(4)] The log normal sheaf to $f$, $\sN^{log}$, is torsion-free.
	\item[(5)] If $|\gamma|\ge 3$, then $f(C)$ is smooth along $E$ and the normal sheaf to $f$ is torsion-free.  In particular, $f(C)$ has no cuspidal singularities.
	\item[(6)] If $|\gamma|\ge 4$, then $f(C)$ has only nodal singularities.
\end{enumerate}
\end{thm}

%{\it Remark.}  If $|\gamma|=2$, then $f(C)$ can be singular along $E$.  For example, if $E$ is given by the equation $y^2 = x^3 - x$, then the morphism $t\mapsto (t^2,t^3)$ defines a cuspical cubic with an order $2$ contact at the origin and an order $7$ contact at infinity.  This lives in a one-dimensional component of $\eK_3^{e_2+e_7}$ and the general curve in this component has a cusp along $E$.  The point is that the torsion section of the normal sheaf to $f$ which would smooth the cusp is not tangent to $E$, so it does not come from a section of the log normal sheaf.  Hence any deformation preserving the contacts will also preserve the cusp.

\begin{proof}  Let $B\subseteq V$ be a representable, smooth, dense open substack.  Such a substack exists because the general stable map in $V$ has no automorphisms.  Let $b\in B$ be a general point.  Let $\pi:C\to B$ be the coarse moduli space of the universal twisted curve over $B$, and let $f:C\to\bP^2$ be the morphism induced by the universal morphism into $\pe$.  Let $f_b:C_b\to\bP^2$ be the restriction to the fiber over $b$.

It follows from Proposition \ref{contact_type} that $f_b^*E$ is supported precisely at the marked points with multiplicities determined according to $\gamma$.  Since the dimension of $V$ is greater than or equal to the expected dimension, which is $|\gamma|-1$, we conclude that
$$|\gamma|-1\le \dim T_bB.$$
Now we prove the reverse inequality.  From Lemma \ref{inf_def} we know that $T_bB$ is a subspace of $H^0(C_b, \sN_b^{log})$, and by Lemma \ref{Cap_Harris_1} we know that it contains no torsion sections.  Since we assumed $|\gamma|\ge 2$, it follows that $\sN_b^{log}/(\sN_b^{log})_{\tors}$ has nonnegative degree on $C_b\cong\bP^1$.  Therefore, $h^1(C_b,\sN_b^{log})=0$ and it follows that $$h^0(C_b,\sN_b^{log})=3d-1-\sum_i (i-1)\gamma_i = |\gamma|-1.$$

We conclude that $\dim T_bB = \dim B = |\gamma|-1$, that $T_bB=H^0(C_b, \sN_b^{log})$, and that $\sN_b^{log}$ is torsion-free.  From this, statements 1, 2, and 4 of the theorem follow.  Now we show that every evaluation map at a twisted marking is surjective.  If not, then we are in the situation of Lemma \ref{Cap_Harris_2} with $f(Q)$ being a single point.  It follows that every section of $\sN_b^{log}$ vanishes to order at least $m$ in $\sN_b/(\sN_b)_{\tors}$, where $m$ is the contact type at this marking.  But this contradicts Lemma \ref{inf_def}, which shows that some section vanishes to order at most $m-1$, at least if $\sN_b^{log}$ has nontrivial sections.  Since $|\gamma|\ge 2$, such a section exists.

Now we show that $\sN_b$ is torsion-free if $|\gamma|\ge 3$.  Since $\sN_b^{log}$ has no torsion, $\sN_b$ can only have torsion at a point $q$ mapping into $E$.  At such a point, $df_b$  vanishes to some order $\ell-1$ with $\ell>1$.  Applying Lemma \ref{Cap_Harris_2}, and the remark following it, we see that either every section of $\sN_b^{log}$ vanishes to order $m-\ell$ in $\sN_b^{log}$ or else every section vanishes to order at least $m$.  The latter cannot happen by the argument of the previous paragraph, and the former could only happen if $h^0(\sN_b^{log})\le 1$.  But this contradicts the assumption $|\gamma|\ge 3$.

Now we show that $f(C)$ is smooth along $E$ if $|\gamma|\ge 3$.  Since $\sN_b$ is torsion-free, a singularity could only occur if two points of $C$ map to the same point of $E$.  But then $\sN_b^{log}$ would have a section vanishing at one of the points and not the other, and the corresponding first order deformation would separate the points.  Since $B$ is smooth at $b$, this deformation extends to a one-parameter family, which contradicts the generality of $b$.

A similar argument can be applied to singularities away from $E$.  If $|\gamma|\ge 4$ and there are three points mapping to the same point of $\bP^2$, then $\sN_b^{log}$ has a section vanishing at two of the points and not the third.  The corresponding deformation moves the third point away from the first two.  This can be applied repeatedly until the image curve has only nodal singularities.
\end{proof}

\subsection{Multiple covers}

We now consider maps whose source curve is smooth but which do not map birationally onto their image.  We need to show that the operation which sends a multiple cover to the normalization of its image defines a morphism, at least on a dense open set.  This allows us to handle irreducible components of the space of stable maps whose general element is a multiple cover if we know enough about the irreducible component corresponding to the image of the multiple cover.

Let $\gamma$ be an infinite sequence with $I\gamma=3d$ for some positive integer $d$.  Let $V\subseteq \eK_d^{\gamma}$ be an irreducible component such that the general point of $V$ corresponds to a map from a smooth source curve which is an $e$ to $1$ cover of its image.  Let $S\to V$ be an \'{e}tale morphism from a smooth integral scheme $S$ such that the image contains only maps from smooth curves which are $e$ to $1$ covers of their image.  This gives us a twisted stable map over $S$.
$$\xymatrix{\Sigma_i \ar@{^{(}->}[r] \ar[d] & \fC \ar[r] \ar[d] \ar@/_/[dd] & \pe \ar[d] \\
\sigma_i \ar@{^{(}->}[r] \ar[dr]_{\cong} & C \ar[d] \ar[r] & \bP^2 \\
 & S,}$$
We now construct an ``image'' stable map over a dense open subvariety of $S$.  Let $\overline{C}$ be the image of $C$ in $S\times\bP^2$, with the reduced induced structure, and let $\overline{C}^{\nu}$ be its normalization.  Note that $\overline{C}$ is flat over $S$, since the fibers are plane curves of degree $d/e$.  It follows that over a dense open subvariety $U\subseteq S$, the fibers of $\overline{C}$ are birational to those of $\overline{C}^{\nu}$.  Moreover, since normal varieties are nonsingular in codimension $1$, there is a dense open subvariety $U'\subseteq U$ such that the pullback of $\overline{C}^{\nu}$ is smooth.  By generic smoothness, there is a $U''\subseteq U'$ over which the fibers are smooth.  Now replace $S$ with $U''$.  We have now constructed the diagram
$$\xymatrix{\overline{C}^{\nu} \ar[r]^f \ar[d]_\pi & \bP^2 \\
S.}$$

The family $\overline{C}^{\nu}$ is flat over $S$ since the fibers are all isomorphic to $\bP^1$.  By the universal property of normalization, we have a factorization $C\to\overline{C}^{\nu}\to\bP^2$.  Let $\overline{\sigma}_1,\ldots, \overline{\sigma}_m$ be the irreducible components of $f^*E$.  These define sections of $\pi$ which are disjoint since the original sections $\sigma_i$ are disjoint.  Let $\varrho_i$ be the multiplicity of $f^*E$ along $\overline{\sigma}_i$.  By Proposition \ref{contact_type}, there is a unique twisted stable map over $S$ of contact type $\varrho_1,\ldots, \varrho_m$ inducing $f:\overline{C}^{\nu}\to\bP^2$.  This defines a morphism $S\to\eK_{d/e}(\vrho)$.

The above construction is used in section 3 to show that components consisting of multiple covers cannot contribute to Gromov-Witten invariants, at least when $|\gamma|\ge 3$.  For now we use it to prove the following.

\begin{thm}
\label{multiple_covers_2}
Let $\varrho_1,\varrho_2,d$ be positive integers such that $\varrho_1+\varrho_2=3d$.  Let $V\subseteq\eK_d(\varrho_1,\varrho_2)$ be an irreducible component whose general source curve is smooth.  Then $V$ is one-dimensional, and $\eK_d(\varrho_1,\varrho_2)$ is generically reduced along $V$.
\end{thm}

\begin{proof}
Let $f:C\to\bP^2$ be the morphism corresponding to a general point of $V$, and suppose it has degree $e$ onto its image.  Since there are only two marked points and every point in $f^*E$ is marked, $f$ must be totally ramified at the marked points.  The above construction defines an \'{e}tale morphism $S\to V$ and a generically injective morphism $S\to \eK_{d/e}(\varrho_1/e,\varrho_2/e)$ (since the ramification points cannot move without moving the image curve).  By Theorem \ref{gen_smooth}, the latter stack is one-dimensional.  Since the expected dimension of $V$ is $1$, it follows that $V$ is also one-dimensional.

By Remark \ref{no_torsion}, $\sN_f^{log}$ has no torsion at the two marked points and from Theorem \ref{gen_smooth}, it follows that $\sN_f^{log}$ has no torsion away from the marked points.  Since $\sN_f^{log}$ has degree $0$ in any case, it must be trivial, and Lemma \ref{inf_def} implies that $\eK_d(\varrho_1,\varrho_2)$ is generically reduced along $V$.
\end{proof}

\subsection{Maps from nodal twisted curves}
\label{sec:nodal}

We review some basic facts about twisted nodes which are used later in the paper.  One application is to show that components of the space of twisted stable maps for which the general source curve has nodes do not contribute to any invariants involving at least two contacts.  We also have a delicate argument in section \ref{sec:large_contact} to deal with points of large contact type on degree $0$ components of the source curve.

First recall how a twisted node can be obtained from its coarse moduli space \cite[\S 3.5]{Abr}.  For a separating node, it is sufficient to use $r$-th root constructions along the two subcurves which the node separates.  While this can be done on the curve itself, it may be easier to visualize it in a one paramater family which smooths the node.  Suppose we have a smooth curve $S$ and a smooth surface $C$ mapping to $S$ whose general fiber is a smooth curve and whose special fiber $C_0$ has a single nodal singularity which separates $C_0$ into two irreducible components.
$$\xymatrix{
C_0 \ar@{^{(}->}[r] \ar[d] & C \ar[d] \\
s_0 \ar@{}[r]|{\in} & S}$$
To make this node be twisted to order $r$, one would apply the $r$-th construction to the two components of $C_0$ (one at a time) and also to the point $s_0$ inside of $S$.  This creates a family $\fC\to\fS$ whose general fiber is the same as before.  The special fiber of the new family (over a morphism $\spec\bC\to\fS$ which maps into $s_0^{1/r}$) is isomorphic to $C_0$ except at the node.  Near the node, the local picture is the same as that of $$[(\spec\bC[x,y]/(xy))/\mu_r],$$
where $\mu_r$ acts on $x$ and $y$ with opposite weights (the \emph{balanced} condition):  $\zeta\cdot(x,y) = (\zeta^{-1}x,\zeta y).$

Now we discuss gluing morphisms, which glue together two marked points on a pair of stable maps.  This requires a precise treatment of evaluation morphisms.  In order to compute Gromov-Witten invariants, it is sufficient to use evaluation morphisms which map to the coarse moduli space.  In other words, sending the stable map in diagram \ref{tw_st_map} to the morphism $S\to\sigma_i\to\bP^2$ defines the evaluation map at the $i$-th marked point.  Note that we cannot map to $\pe$, because the morphism $\Sigma_i\to S$ need not have a section.

The more precise evaluation maps, which were first described in \cite{AGVold}, send diagram \ref{tw_st_map} to the subdiagram
$$\xymatrix{\Sigma_i \ar[r] \ar[d] & \pe \\
S}$$
which is an object of the stack of gerbes in $\pe$ \cite[3.3]{AGV}.  This is isomorphic to a rigidification of the inertia stack, $\inertia(\pe)$ \cite[\S 3.4]{AGV}.  The coarse moduli space of $\inertia(\pe)$ is a disjoint union of $\bP^2$ with $r-1$ copies of $D$.  By suitably numbering these components, the evaluation morphism determines the contact type of the marked point through the connected component to which it maps.

Using these evaluation maps, for contact types $\vrho^{(1)}$ and $\vrho^{(2)}$ and an integer $1\leq k\leq r-1$, there is a gluing morphism \cite[\S 5.2]{AGV}
$$eK_{d_1}(\vrho^{(1)},k)\times_{\inertia(\pe)} \eK_{d_2}(\vrho^{(2)},r-k) \to \eK_{d_1+d_2}(\vrho^{(1)},\vrho^{(2)}).$$
Here it is important that one of the evaluation maps is composed with the involution $\inertia(\pe)\to\inertia(\pe)$ which sends a contact type $k$ to its complement $r-k$ for $k>0$.  This is due to the balanced condition at twisted nodes.

Finally, we mention that the image of the gluing morphism has codimension at most $1$ everywhere.  By deformation theory, it can be shown that the fiber of the normal sheaf at a point is at most $1$-dimensional, generated by a first order deformation of the map which infinitesimally smooths the node \cite[p.100]{HM}, \cite[\S 4.6]{AGVold}.  It is useful to note that a first order deformation smoothing the twisted node does not smooth the node in its coarse moduli space.  To see this, recall the construction from the second paragraph of this subsection.  The deformation space of the twisted node, $\fS$, is ramified over the deformation space of the coarse moduli space.  Therefore, one must deform the twisted node to higher order before the coarse moduli space deforms.

\numberwithin{equation}{section}

\section{Enumerativity}

First we define the numbers we wish to compute.  Let $\alpha$ and $\beta$ be two infinite sequences of integers and $d$ a positive integer such that $I\alpha+I\beta=3d$.  We assume $|\beta|>0$.  For each positive integer $i$, fix $\alpha_i$ general points $r_{ij}\in E$.  Also fix $|\beta|-1$ general points $p_i\in\bP^2$.

\begin{definition}
\label{def_Nd}
Let $N_d(\alpha,\beta)$ be the number of rational degree $d$ curves $C$ passing through each point $p_i$, meeting $E$ in an $i$-th order contact at each point $r_{ij}$, having $\beta_i$ additional $i$-th order contacts with $E$, and meeting $E$ only in points where $C$ is unibranch.
\end{definition}

{\it Remark.}  We allow unibranch singularities so that the Gromov-Witten invariants  compute these numbers when $|\alpha|+|\beta|=2.$  When $|\alpha|+|\beta|>2$, all contacts with $E$  occur at smooth points (for generic initial data), and when $|\alpha|+|\beta|>3$, all such $C$  have at worst nodal singularities.
\smallskip

In this section, we compare $N_d(\alpha,\beta)$ with the Gromov-Witten invariants $I_d^r(\alpha,\beta)$.  Fix infinite sequences of integers $\alpha$ and $\beta$ and let $d$ and $r$ be positive integers with $r>3d$.  Assume that $k:=(I\alpha+I\beta-3d)/r$ is an integer and that $|\beta|>k$.  Let $\gamma=(|\beta|-k-1)e_0+\alpha+\beta$.  For $1\le j\le \alpha_i$, let $e_{ij}:\eK_d^{\gamma}\to E$ be the evaluation map at the $j$-th point which has contact type $i$.  For $0\le i\le r-1$, let $p_i$ be the class of a point in the $i$-th component of $\inertia(\pe)$ (which is $\pe$ if $i=0$ and a gerbe over $E$ if $i>0$).

\begin{definition}
\label{GWI_def}
$$I^r_d(\alpha,\beta) = \int_{[\eK_d^{\gamma}]^{vir}} \prod_{j=1}^{|\beta|-1} e_{0j}^*p_0\,\, \prod_{i=1}^{3d}\prod_{j=1}^{\alpha_i} e_{ij}^*p_i$$
\end{definition}

We are interested in the invariants when $k=0$, but other invariants arise in our WDVV recursion.  Let $$ev:\eK_d^{\gamma}\to (\bP^2)^{\gamma_0}\times E^{|\gamma|-\gamma_0}$$ be the product of all evaluation maps.  In order for an irreducible component $V\subseteq \eK_d^{\gamma}$ to be seen by the integral (\ref{GWI_def}), it is necessary that $ev(V)$ have dimension large enough to support the pushforward of the virtual fundamental class.  In other words, it is necessary that $\dim(ev(V))\ge\edim(V)$, where $\edim(V)$ is the expected dimension which $V$ inherits from $\eK_d^{\gamma}$.  So our strategy for demonstrating that the Gromov-Witten invariants are enumerative is to first show that this fails for multiple cover components, and then to show that it fails for components in which the general source curve is singular.

\begin{lemma}
\label{smoothcurve}
Assume $I\gamma=3d$ and let $V\subseteq\eK_d^{\gamma}$ be an irreducible component for which the general map $f:C\to\bP^2$ has smooth source curve $C$.  Then $\dim(ev(V))\le \edim(V)$ with equality only if one of the following holds.
\begin{itemize}
	\item $C\to f(C)$ is birational
	\item $|\gamma|-\gamma_0=2$ and $C\to f(C)$ is a cyclic cover, ramified at the twisted marked points
	\item $|\gamma|-\gamma_0=1$
\end{itemize}
\end{lemma}

\begin{proof}
If the general map is an $e$ to $1$ cover of its image, then the construction of section 2.3 defines an \'{e}tale morphism $\phi:S\to V$ and a morphism $\psi:S\to \eK_{d/e}^{\delta}$, where $\delta$ is determined by the general image curve.  Let $W=\overline{\psi(S)}.$  Each twisted marked point on the image curve corresponds to a collection of twisted marked points on $C$, and under this correspondence, $\phi$ and $\psi$ are compatible with the evaluation maps.  It follows that $\dim(ev(V))\le\dim(ev(W)).$  We also have $\edim(W)\le\edim(V)$, since $C$ has at least as many marked points as its image.  We claim that $\dim(ev(W))\le\edim(W).$  If $|\gamma|-\gamma_0\ge 2$, then Theorem \ref{gen_smooth} implies that $\dim(W)=\edim(W)$, from which the inequality follows.  If the image curve has a single twisted marked point, then it can only have a single contact with $E$, which must occur at a $3{d/e}$ torsion point.  
Moreover, the image of the universal curve over $W$ in $\bP^2$ is a Zariski-closed subset containing a finite number of points of $E$, so it must be one-dimensional.  It follows that $\dim(ev(W))=\gamma_0=\edim(W)$ in this case.

We conclude that $\dim(ev(V))\le\edim(V).$  If they are equal, then each of the three inequalities above is an equality, and in particular $\edim(W)=\edim(V)$.  This implies that both $C$ and its image have the same number of twisted markings.  If $e>1$, this is only possible if there are at most two marked points.  Moreover, when there are precisely two marked points, $C\to f(C)$ must be a cyclic cover, since it has to be totally ramified at the two twisted markings.
\end{proof}

\begin{lemma}
\label{nodalcurve}
Assume $I\gamma=3d$ and $r>3d$.  If $V\subseteq\eK_d^{\gamma}$ is an irreducible component, then $\dim(ev(V))\le\edim(V)$.  Moreover, if $\dim(ev(V))=\edim(V)$ and $|\gamma|-\gamma_0\ge 2$, then the general map in $V$ has irreducible source curve, $V$ has the expected dimension, and $\eK_d^{\gamma}$ is generically reduced along $V$.
\end{lemma}

\begin{proof}
Recall that we implicitly assume $r>3d$.  The proof is by induction on the number of irreducible components of the general source curve in $V$.  Note that untwisted markings are easy to deal with, because the universal curve over $\eK_d^{\gamma}$ is isomorphic to $\eK_d^{\gamma+e_0}$.  Therefore, $\dim(ev(V))$ increases by one whenever an untwisted marking is added to a component mapping with positive degree.

Let $f:C\to\bP^2$ be a general map in $V$.  If $C$ is irreducible, then $\dim(ev(V))\le\edim(V)$ by Lemma \ref{smoothcurve}.  To prove the last statement, we treat two cases.  If $|\gamma|-\gamma_0=2$, then the statement follows from Theorem \ref{multiple_covers_2}.  If $|\gamma|-\gamma_0>2$, then Lemma \ref{smoothcurve} implies that $C\to f(C)$ is birational, so the statement follows from Theorem \ref{gen_smooth}.

Now suppose that $C$ has a node $\eta$ between two components which map with positive degree.  This node must be untwisted, because otherwise the contact types would sum to $r$, which is impossible since $r>3d$.  Let $V_1$ and $V_2$ be the irreducible components of twisted stable map spaces which contain the maps obtained by separating $C$ at $\eta$.  Then we have the following diagram.
$$\xymatrix{
 & & ev(V_1) \times ev(V_2) \ar[d] \ar@{^{(}->}[r] & (\bP^2)^{\gamma_0+2}\times E^{|\gamma|-\gamma_0} \ar[d] \\
V \ar[r] & V_1\times V_2 \ar[r] \ar[ru] & ev(V) \ar@{^{(}->}[r] & (\bP^2)^{\gamma_0}\times E^{|\gamma|-\gamma_0}}$$
Since the fibers of the vertical arrows are $2$-dimensional, we have by induction $$\dim(ev(V)) = \dim(ev(V_1))+\dim(ev(V_2)) -2 \le \edim(V_1)+\edim(V_2) -2 = \edim(V) -1.$$
This proves the lemma for $V$.

Now suppose that $C$ has a node between two components which map with degree $0$.  Then the node could be smoothed, because $\eK_{0,n}(B\mu_r)$ is flat over $\overline{\eM}_{0,n}$ \cite[3.0.5]{ACV}.  So for a general $f:C\to\bP^2$, in $V$, such a node does not exist.  Finally, suppose that $C$ has a degree $0$ component $C_0$ containing $s$ marked points, and let $C_1,\ldots, C_m$ be the connected components of $\overline{C\setminus C_0}$.  Let $V_i$ be the irreducible component containing the map $C_i\to\bP^2$.  We first treat the case where $C_0$ maps away from $E$.  Now by a similar argument as above, this time using the fact that each of the marked points where $C_i$ connects to $C_0$ are constrained to map to a single point, we have
$$\dim(ev(V))\le \sum_{i=1}^m \dim(ev(V_i)) - m + 1.$$
By induction, this implies that
$$\dim(ev(V))\le \sum_{i=1}^m \edim(V_i) - m + 1 = \edim(V) - s - m + 2.$$
Since $s+m\ge 3$ by stability, the Lemma is proven for $V$.

Now we treat the case where $C_0$ maps to $E$, and for simplicity we assume that $\gamma_0=0$.  In this case, the point where $C_i$ meets $C_0$ is twisted, hence cannot move around on $C_i$, so we have the weaker inequality
$$\dim(ev(V))\le\sum_{i=1}^m\dim(ev(V_i)),$$
which by induction yields
$$\dim(ev(V))\le\sum_{i=1}^m\edim(V_i)=\edim(V)-s+1.$$
We claim that $s\ge 1$.  If not, then by Lemma \ref{contact_deg0}, the contact types at the nodes of $C_0$ would have to sum to a multiple of $r$.  By the balanced condition, the same would be true of the orders of contact between $C_i$ and $E$ at the nodes.  But this contradicts $r>3d$.

Therefore, $\dim(ev(V))\le\edim(V)$ with equality if and only if the following three conditions hold.
\begin{itemize}
	\item $s=1$
	\item $\dim(ev(V))=\sum_{i=1}^m\dim(ev(V_i))$
	\item $\dim(ev(V_i))=\edim(V_i)$
\end{itemize}
We claim that under these assumptions, each $C_i$ has a single marked point where it meets $C_0$.  If not, then $C_i$ is irreducible by induction.  If $C_i$ maps birationally onto its image, then Theorem \ref{gen_smooth} implies that the evaluation map $V_i\to E$ at the node with $C_0$ is surjective.  If it is a multiple cover of its image, then Lemma \ref{smoothcurve} implies that there is a bijection between the twisted markings of $C_i$ and those of its image, which again implies surjectivity of this evaluation map.  But this means that the node imposes a condition on $\times_i ev(V_i)$ ($m\ge 2$ by stability), which contradicts the second condition above.  Therefore, the curve $C$ has a single marked point, which lies on $C_0$.  This verifies the lemma for $V$.
\end{proof}

We are ready to prove that the invariants are enumerative, at least when $|\alpha|+|\beta|\ge 2$.

\begin{thm}
\label{enumerative}
If $3d=I\alpha+I\beta$, $r>3d$, and $|\alpha|+|\beta|\ge 3$, then $$I_d^r(\alpha,\beta)=\beta!N_d(\alpha,\beta).$$  If $|\alpha|+|\beta|=2$, then there are multiple cover contributions which are accounted for by the formulas
$$I_d^r(e_k,e_{\ell})=\sum_{s|\gcd(d,k,\ell)} \frac{1}{s}N_{d/s}(e_{k/s},e_{\ell/s}),$$
and 
$$\frac{I_d^r(0,e_k+e_{\ell})}{(e_k+e_{\ell})!}=\sum_{s|\gcd(d,k,\ell)} N_{d/s}(0,e_{k/s}+e_{\ell/s}).$$
\end{thm}

\begin{proof}
Let $\gamma=(|\beta|-1)e_0+\alpha+\beta$.  First we treat the case $|\alpha|+|\beta|\ge 3$.  By Lemmas \ref{smoothcurve} and \ref{nodalcurve}, the only irreducible components $V\subseteq\eK_d^{\gamma}$ which contribute to $I_d(\alpha,\beta)$ are those for which the general curve is smooth and maps birationally onto its image.  Moreover, these components have the expected dimension and $\eK_d^{\gamma}$ is generically reduced along them.  Therefore, $I_d(\alpha,\beta)$ is equal to the sum over such components $V$ of the degrees of the evaluation morphisms $$V\to (\bP^2)^{\gamma_0}\times E^{|\alpha|}.$$  By generic smoothness, this equals the number of preimages of a general point, which recovers Definition \ref{def_Nd} in light of Theorem \ref{gen_smooth}.  The factor of $\beta!$ is due to the fact that the twisted marked points having contact type $k$ which are not constrained to map to a given point can be permuted.

If $|\gamma|=2$, then the only difference is that there can be multiple cover contributions.  We can still apply generic smoothness in the same context, but now the degree has a factor of $1/s$ if the general map in $V$ is an $s$ to $1$ cover of its image.  This factor comes from the automorphism group of the cover, which has order $s$ since it is cyclic by Lemma \ref{smoothcurve}.  In the last formula, the factor of $1/s$ disappears because there is an untwisted marking on the curve.
\end{proof}

\section{Application of WDVV}

We wish to use the WDVV equations to compute relations among the twisted Gromov-Witten invariants $I_d^r(\alpha,\beta)$.  For a precise statment and proof of WDVV in the context of Deligne-Mumford target stacks, see \cite[\S 6.2]{AGV}.  Fix integers $r$ and $d$ such that $r>6d$.  For $1\le i\le r-1$, let $f_i$ denote the fundamental class of the $i$-th component of the inertia stack of $\pe$, and let $p_i$ denote the class of a point on this component.  Let $h$ denote the hyperplane class in $\bP^2$.  Fix an integer $k$ with $1\le k\le 3d-1$, and choose sequences $\alpha,\beta$ so that $I\alpha+I\beta+k+1=3d$.

We implement WDVV using the sequence of classes $f_k,f_1,h,h$.  Let $n=|\alpha|+2|\beta|$ and let $\delta=(\delta_1,\ldots,\delta_n)$ be an $n$-tuple of Chow classes consisting of
\begin{itemize}
	\item $\alpha_i$ copies of $p_i$ for $1\le i\le r-1$,
	\item $\beta_i$ copies of $f_i$ for $1\le i\le r-1$, and
	\item $|\beta|$ copies of $p$.
\end{itemize}
We use the notation of \cite[\S 6.1]{AGV}, except that we let $$\langle\gamma_1,\ldots,\gamma_m\rangle_d = \int_{[\eK_{0,m}(\pe,d)]^{vir}} \prod_{i=1}^m e_i^*\gamma_i.$$  By a standard argument, \cite[6.2.1]{AGV} becomes the statement that for any four classes $\gamma_1,\ldots,\gamma_4$, the quantity
\begin{gather*}
\sum_{d_1+d_2=d}\sum_{A\sqcup B=\{1,\ldots,n\}} \left(\sum_{i=0}^2 \langle\gamma_1,\gamma_2,\delta_A,h^i\rangle_{d_1} \langle\gamma_3,\gamma_4,\delta_B,h^{2-i}\rangle_{d_2}+\right. \\
\left. r\sum_{i=1}^{r-1} \langle\gamma_1,\gamma_2,\delta_A,f_i\rangle_{d_1} \langle\gamma_3,\gamma_4,\delta_B,p_i\rangle_{d_2} + r\sum_{i=1}^{r-1} \langle\gamma_1,\gamma_2,\delta_A,p_i\rangle_{d_1} \langle\gamma_3,\gamma_4,\delta_B,f_i\rangle_{d_2}\right)
\end{gather*}
is invariant under interchanging $\gamma_2$ and $\gamma_3$.

Applying this to the four classes $f_k,f_1,h,h$ and using the $n$-tuple $\delta$ defined above, we obtain the following equation.  We use a bare summation to indicate a sum over $d_1+d_2=d$, $\alpha^1+\alpha^2=\alpha$, and $\beta^1+\beta^2=\beta$, with $0<d_1<d$.

%\numberwithin{figure}{section}
%\begin{figure}[htb]
\begin{gather}
d^2I_d(\alpha,\beta+e_{k+1})= 3dI_d(\alpha+e_1,\beta+e_k)+3dI_d(\alpha+e_k,\beta+e_1)-I_d(\alpha,\beta+e_1+e_k) \notag \\
+
\sum I_{d_1}(\alpha^1,\beta^1+e_1)I_{d_2}(\alpha^2,\beta^2+e_k) d_1^2d_2^2\ch{\alpha}{\alpha^1}\ch{\beta}{\beta^1}\ch{|\beta|}{|\beta^1|} \label{WDVV_equation} \\
-
\sum I_{d_1}(\alpha^1,\beta^1+e_1+e_k)I_{d_2}(\alpha^2,\beta^2) d_1d_2^3\ch{\alpha}{\alpha^1}\ch{\beta}{\beta^1}\ch{|\beta|}{|\beta^1|+1} \notag \\
+
r\sum_{\ell=1}^{r-1}\sum I_{d_1}(\alpha^1+e_{\ell},\beta^1+e_1)I_{d_2}(\alpha^2,\beta^2+e_k+e_{r-\ell}) d_1d_2\ch{\alpha}{\alpha^1}\ch{\beta}{\beta^1}\cdot
\left\{\begin{array}{cc} \ch{\vert\beta\vert}{\vert\beta^1\vert} & \ell\le 3d \notag \\
\ch{\vert\beta\vert}{\vert\beta^2\vert+1} & \ell>3d \end{array}\right. \notag \\
+
r\sum_{\ell=1}^{r-1}\sum I_{d_1}(\alpha^1,\beta^1+e_1+e_{\ell})I_{d_2}(\alpha^2+e_{r-\ell},\beta^2+e_k) d_1d_2\ch{\alpha}{\alpha^1}\ch{\beta}{\beta^1}\cdot
\left\{\begin{array}{cc} \ch{|\beta|}{|\beta^1|+1} & \ell\le 3d \notag \\
\ch{|\beta|}{|\beta^2|} & \ell>3d \end{array}\right. \notag \\
-
r\sum_{\ell=1}^{r-1}\sum I_{d_1}(\alpha^1+e_{\ell},\beta^1+e_1+e_k)I_{d_2}(\alpha^2,\beta^2+e_{r-\ell}) d_2^2\ch{\alpha}{\alpha^1}\ch{\beta}{\beta^1}\cdot
\left\{\begin{array}{cc} \ch{|\beta|}{|\beta^1|+1} & \ell\le 3d \notag \\
\ch{|\beta|}{|\beta^2|} & \ell>3d \end{array}\right. \notag \\
-
r\sum_{\ell=1}^{r-1}\sum I_{d_1}(\alpha^1,\beta^1+e_1+e_k+e_{\ell})I_{d_2}(\alpha^2+e_{r-\ell},\beta^2) d_2^2\ch{\alpha}{\alpha^1}\ch{\beta}{\beta^1}\cdot
\left\{\begin{array}{cc} \ch{|\beta|}{|\beta^1|+2} & \ell\le 3d \notag \\ \ch{|\beta|}{|\beta^2|-1} & \ell>3d \end{array}\right.
\end{gather}
%\caption{A WDVV equation for $\pe$.}
%\end{figure}

At the end of the paper, we show how this can be combined with our other formulas to compute the numbers $N_d(\alpha,\beta)$.

\section{Dealing with a large contact type}

\label{sec:large_contact}

In the recursion above, the last four sums involve a large contact type, which must lie on a component mapping with degree $0$.  By large contact type, we simply mean a contact type which is larger than $3d$, the largest contact type that can appear on a component which maps with positive degree at most $d$.  The invariants which involve a single large contact type can be reduced to invariants without large contacts via the following formulas.

\begin{thm}
\label{large_contact_type}
If $r>I\alpha+I\beta=3d+\ell$, $3d>\ell$, and $\alpha_i=\beta_i=0$ for $i>3d$, then
\begin{align*}
rI_d(\alpha+e_{r-\ell},\beta) &= \sum_{k=\ell+1}^{3d} (k-\ell)\beta_k I_d(\alpha+e_{k-\ell},\beta-e_k), \\
rI_d(\alpha,\beta+e_{r-\ell}) &=\sum_{k=\ell+1}^{3d}(k-\ell)\beta_k
I_d(\alpha,\beta+e_{k-\ell}-e_k) +
 \sum_{k=\ell+1}^{3d} (k-\ell)\alpha_k I_d(\alpha+e_{k-\ell}-e_k,\beta).
\end{align*}
\end{thm}

\begin{proof}
Choose an integer $k > \ell$ with $k\le 3d$, and let $\gamma = \alpha + \beta$.  Note that under the hypotheses $I\gamma-\ell=3d$ and $\gamma_i=0$ for $i>3d$, it is necessary that $|\gamma|\ge 2$.  We consider the following diagram.

$$\xymatrix{
& {\displaystyle \coprod^{\gamma_k} \eK_d^{\gamma+e_{k-\ell}-e_k}\times_{\inertia(\pe)} \eK_0(r+\ell-k, k, r-\ell)} \ar[dl]_{F_k} \ar[dr]^{q_k} & \\
{\displaystyle\eK_d^{\gamma+e_{r-\ell}}} & & {\displaystyle \coprod^{\gamma_k} \eK_d^{\gamma+e_{k-\ell}-e_k}}}$$
\medskip

The morphisms to $\inertia(\pe)$ in the fiber product are evaluation maps at the marked points corresponding to the extra $k-\ell$ contact type and the $r+\ell-k$ contact type, respectively (one of these is composed with the involution of $\inertia(\pe)$ which interchanges contact types summing to $r$).  The morphism $F_k$ is the gluing morphism associated to these two points \cite[5.2.1.1]{AGV}, and $q_k$ is the first projection.  We need to take a disjoint union of $\gamma_k$ copies of the fiber product, because there are $\gamma_k$ points with contact type $k$ which could be on the degree $0$ component.  Note that these morphisms are compatible with the evaluation maps.

We first make the followings claims.
\begin{enumerate}
	\item The morphism $q_k$ is an \'{e}tale gerbe banded by a cyclic group of order $r/\gcd(r,k-\ell)$.
	\item For any irreducible component $$V\subseteq \eK_d^{\gamma+e_{k-\ell}-e_k}\times_{\inertia(\pe)} \eK_0(r+\ell-k, k, r-\ell),$$
$(F_k)_*[V] = [F_k(V)].$  Here $V$ and $F_k(V)$ have the reduced induced structure.
	\item If $q_k(V)$ contributes to Gromov-Witten invariants, then $F_k(V)$ is an irreducible component.
	\item If an irreducible component $W\subseteq \eK_d^{\gamma+e_{r-\ell}}$ contributes to Gromov-Witten invariants, then $W=F_k(V)$ for some integer $k$ and irreducible component $V$.
	\item If either $q_k(V)$ or $F_k(V)$ contributes to Gromov-Witten invariants, then
	\begin{enumerate}
		\item $q_k(V)$ and $F_k(V)$ have the expected dimension,
		\item $\eK_d^{\gamma+e_{k-\ell}-e_k}$ is generically reduced along $q_k(V)$, and
		\item $\eK_d^{\gamma+e_{r-\ell}}$ has multiplicity $(k-\ell)/\gcd(r,k-\ell)$ along $F_k(V)$.
	\end{enumerate}
\end{enumerate}

\noindent CLAIM 1:  We will construct an isomorphism $\eK_0(r+\ell-k,k,r-\ell)\to E^{1/r}$ such that the evaluation map $E^{1/r}\to\inertia(\pe)$ at the first marking is the same as the rigidification map which ``quotients out'' the cyclic group which stabilizes this marking.  Since this stabilizer group has order $r/\gcd(r,k-\ell)$, the claim follows by taking a base change with the other evaluation map $\eK_d^{\gamma+e_{k-\ell}-e_k}\to \inertia(\pe)$.  Before constructing this isomorphism, note that Lemma \ref{contact_deg0} gives us a bijection between the $\bC$-valued points of both stacks.  Moreover, all of these points have stabilizer group $\mu_r$.  Since both stacks are smooth, it suffices to show that there is a morphism from one to the other inducing this bijection which is an isomorphism on stabilizer groups.

To construct the morphism, start with a trivial $\bP^1$ bundle $\bP^1\times E$ and apply root constructions to three constant sections (call them $\sigma_i$, and let $\sigma$ be a general constant section), to obtain a family of three-pointed twisted curves over $E$.  The orders of the roots are determined by the three contact types.  To map this family to $E^{1/r}$, we need an $r$th root of $p_2^*\sO_E(E)$.  We have an $r$th root of the trivial bundle, given by
$$\sO\left(\frac{r+\ell-k}{r}\sigma_1 + \frac{k}{r}\sigma_2 + \frac{r-\ell}{r}\sigma_3 - 2\sigma\right),$$ which has the correct contact types at the marked points.  If we pull everything back to the gerbe $E^{1/r}\to E$, then we can tensor this line bundle with $p_2^*\sO_{E^{1/r}}(E^{1/r})$ to obtain the required $r$th root.  This gives us a family
$$\xymatrix{\sP \ar[r] \ar[d] & E^{1/r} \ar@{^{(}->}[r] & \pe \\
E^{1/r}},$$
which defines the required morphism $E^{1/r}\to \eK_0(r+\ell-k,k,r-\ell)$.
\medskip

\noindent CLAIM 2:  It suffices to show that $F_k$ defines an injection between $\bC$-valued objects and identifies their automorphism groups.  This follows from the results of \cite{AGV}, particularly Proposition A.0.2.
\medskip

\noindent CLAIM 3:  First we remark that $\coprod F_k$ is not surjective, because the large contact type point can be contained on a degree zero component having more than one other marked point and meeting more than one other component.  To get a surjective map, we would have to take a disjoint union over more than just the marked points.  If we did this, then the images of the irreducible components would be irreducible closed subsets of $\eK_d^{\gamma+e_{r-\ell}}$ whose union is the whole space.  An irreducible component of $\eK_d^{\gamma+e_{r-\ell}}$ is therefore a maximal element of this collection.

If $q_k(V)$ contributes to Gromov-Witten invariants, then Lemma \ref{nodalcurve} implies that the general source curve in $V$ is irreducible, and it follows that $F_k(V)$ is an irreducible component.
\medskip

\noindent CLAIM 4:  Let $f:C\to\bP^2$ be a general map in $W$ and let $C_0\subseteq C$ be the component containing the $r-\ell$ contact type point.  Suppose this component contains $s$ marked points and is connected to $m$ other components.  As in the proof of Lemma \ref{nodalcurve}, it follows that $$\dim(ev(W))\le \edim(W) -s+2.$$  It is $+2$ instead of $+1$, because the expected dimension formula is altered by the large contact type.  Therefore $s\le 2$, and we also have $s\ge 2$ by the following argument.  We know that $C_0$ contains the $r-\ell$ contact type point, and if all other twisted points were nodal, their contact types would have the form $r-k_i$, with $\sum k_i\le 3d$.  But the contact types must sum to a multiple of $r$, which is impossible since $r>3d+\ell$.  We claim that there is only one component meeting $C_0$.  
If not, then $\dim(ev(V))$ would be reduced by $1$ due to the condition imposed by two nodes mapping to the same point of $E$, which is a nontrivial condition since $|\gamma|\ge 2$ implies that one of the evaluation maps at the node is surjective by Theorem \ref{gen_smooth} and Lemma \ref{smoothcurve}.  It follows that $W=F_k(V)$ for some integer $k$ and irreducible component $V$.
\medskip

\noindent CLAIM 5:  First note that for $q_k(V)$ and $F_k(V)$, the quantities $\dim(\cdot)$, $\edim(\cdot)$, and $\dim(ev(\cdot))$ agree.  Therefore, to prove part (a) of the claim, it suffices to prove that if $\dim(ev(q_k(V)))\ge\edim(q_k(V))$, then $q_k(V)$ has the expected dimension.  But this follows from Lemma \ref{nodalcurve}, since $|\gamma|\ge 2$.  Part (b) of the claim follows from the same Lemma.

Recall from section 2.4 that the normal space to the image of $F_k$ in $\eK_d^{\gamma+e_{r-\ell}}$ consists of first order deformations which smooth the node.  Since $\eK_d^{\gamma+e_{k-\ell}-e_k}$ is generically smooth along $q_k(V)$, the only contribution to the multiplicity of $\eK_d^{\gamma+e_{r-\ell}}$ along $F_k(V)$ comes from the normal space to the image of $F_k$, hence from deformations smoothing the node.

Recall that the deformation space of a twisted node is obtained from the coarse moduli space by applying root constructions, as we discussed in section 2.4.  So we start with the family $xy=t$ (suitably compactified) with two sections $\sigma_1$ and $\sigma_2$ given by $y=a$ and $y=b$.  Let $X$ be the component $x=0$ of the central fiber, let $Y$ be $y=0$, and let $\sigma$ be another section given by $y=c$.  

\numberwithin{figure}{section}
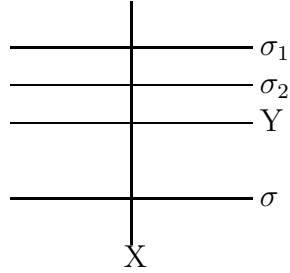
\begin{figure}[htb]
\setlength{\unitlength}{2500sp}%
\begingroup\makeatletter\ifx\SetFigFont\undefined%
\gdef\SetFigFont#1#2#3#4#5{%
  \reset@font\fontsize{#1}{#2pt}%
  \fontfamily{#3}\fontseries{#4}\fontshape{#5}%
  \selectfont}%
\fi\endgroup%
\begin{picture}(2625,2637)(4789,-2986)
\put(7276,-2386){\makebox(0,0)[lb]{\smash{{\SetFigFont{12}{14.4}{\rmdefault}{\mddefault}{\updefault}{$\sigma$}%
}}}}
\thinlines
{\put(4801,-1561){\line( 1, 0){2400}}
}%
{\put(4801,-811){\line( 1, 0){2400}}
}%
{\put(4801,-1186){\line( 1, 0){2400}}
}%
{\put(4801,-2311){\line( 1, 0){2400}}
}%
\put(5926,-2986){\makebox(0,0)[lb]{\smash{{\SetFigFont{12}{14.4}{\rmdefault}{\mddefault}{\updefault}{X}%
}}}}
\put(7276,-1636){\makebox(0,0)[lb]{\smash{{\SetFigFont{12}{14.4}{\rmdefault}{\mddefault}{\updefault}{Y}%
}}}}
\put(7276,-1261){\makebox(0,0)[lb]{\smash{{\SetFigFont{12}{14.4}{\rmdefault}{\mddefault}{\updefault}{$\sigma_2$}%
}}}}
\put(7276,-886){\makebox(0,0)[lb]{\smash{{\SetFigFont{12}{14.4}{\rmdefault}{\mddefault}{\updefault}{$\sigma_1$}%
}}}}
{\put(6001,-361){\line( 0,-1){2400}}
}%
\end{picture}%
\caption{Family of curves with central fiber $X\cup Y$ and sections.}
\end{figure}

Applying suitable root constructions to $\sigma_1$, $\sigma_2$, $X$, and $Y$, we obtain a family of twisted curves over the affine line with a root construction at $t=0$.  For example, along $X$, $Y$, and $t=0$, the $r'$-th root construction is applied, where $r'=r/\gcd(r,k-\ell)$.  Consider the line bundle $$\sL = \sO\left(\frac{k-\ell}{r} X + \frac{k}{r}\sigma_1 + \frac{r-\ell}{r}\sigma_2 - \sigma\right).$$  We want to map this family of curves into $\pe$ in such a way that $\sL$ is the pullback of $\sO_{\pe}(E^{1/r})$.  First, choose a general map $f:\fC\to\pe$ in $F_k(V)$.  We want to map the central fiber of our family according to $f$, and for this $\sL$ must be chosen as above in order to have the correct contact types.  The pullback of the section $s_{E^{1/r}}$ vanishes along $X$.  It is impossible to extend this section to the whole family (which was to be expected, since $V\to F_k(V)$ is bijective), but it can be extended to an infinitesimal neighborhood of the central fiber.  
Along $X$, such an extension amounts to an infinitesimal deformation of the zero section of $\sL|_X^{1/r'}$.  Since this bundle has no sections, the extension must be trivial.  Along $Y$, the section is $\tau^{\delta/r'}$, where $\delta = (k-\ell)/\gcd(r,k-\ell)$, and $\tau$ is the tautological section of $\sO(X^{1/r'}).$  In order for $\tau^{\delta/r'}$ to be zero along the infinitesimal neighborhood of $X$, we must set $t^{\delta/r'}=0$ on the base.

We then obtain a family of twisted curves over $$\left[\left(\spec\bC[t^{1/r'}]/(t^{\delta/r'})\right)/\bmu_{r'}\right].$$
Note that the coarse moduli space of this family is the $t=0$ fiber of the original family.  Therefore, it already maps into $\bP^2$ on the level of coarse moduli spaces (induced from the original map $f$), and it is only necessary to lift this to the stack.  This is essentially done using $\sL$ and $\tau^{\delta/r'}$.  We conclude that we can deform to order $\delta$ in the normal direction to $V$ inside of $F_k(V)$.  This concludes the proof of claim 5.
\medskip

>From these claims, it is not hard to verify the formulas.  The compatibility of the evaluation maps means that all the invariants can be computed in $$\coprod_{k=\ell+1}^{3d} \coprod^{\gamma_k} \eK_d^{\gamma+e_{k-\ell}-e_k}\times_{\inertia} \eK_0(r+\ell-k, k, r-\ell),$$
provided that the correct analogue of the virtual fundamental class is used.  Since $q_k$ is an \'{e}tale gerbe, the invariants on the right hand side of the formulas are computed by integrating over $$q_k^*[\eK_d^{\gamma+e_{k-\ell}-e_k}]^{vir}\cdot r/\gcd(r,k-\ell).$$
Moreover, claim 5 implies that $[\eK_d^{\gamma+e_{k-\ell}-e_k}]^{vir}$ restricts to $[q_k(V)]$ on any component $q_k(V)$ which contributes to Gromov-Witten invariants.  Claims 2 and 5 imply that $$(F_k)_*([V]\cdot r/\gcd(r,k-\ell)) = [F_k(V)]\cdot r/\gcd(r,k-\ell) = [\eK_d^{\gamma+e_{r-\ell}}]|_{F_k(V)}\cdot r/(k-\ell).$$  The formulas now follow by summing over irreducible components, using claims 3 and 4.
\end{proof}

\section{Using the Caporaso-Harris formula}

In this section we prove an analogue of a formula of Caporaso and Harris for rational plane curves having prescribed tangencies with a smooth plane cubic.  Lemma \ref{relations} and Theorem \ref{CH_recursion} can be generalized to smooth plane curves of degree $\ge 3$.

Fix positive integers $r$ and $d$, and an $n$-tuple of integers $\varrho_1,\ldots,\varrho_n$ such that $1\le\varrho_i\le r-1$ and $\sum\varrho_i = 3d$.  Note that all markings in $\eK_d(\vrho)$ are twisted.  Let $A^*(\eK_d(\vrho))$ be the operational Chow ring of $\eK_d(\vrho)$ \cite{Vi}.  We define $N^1(\eK_d(\vrho))$ to be $A^1(\eK_d(\vrho))$ modulo the equivalence relation $a_1\equiv a_2$ if for all $b\in A_1(\eK_d(\vrho))$, $$\int_{\eK_d(\vrho)}a_1\cap b=\int_{\eK_d(\vrho)}a_2\cap b.$$  Let $\eC$ be the universal coarse curve over $\eK_d(\vrho)$, so that we have the diagram
$$\xymatrix{
\eC \ar[r]^f \ar[d]^{\pi} & \bP^2 \\
\eK_d(\vrho),}$$
with $\pi$ representable.

We define classes $h,\chi_1,\ldots,\chi_n\in N^1(\eK_d(\vrho))$ as follows.  Let $h=\pi_*(f^*p)$, where $p\in A^2(\bP^2)$ is the class of a point and $\pi_*$ is the proper, flat pushforward.  In the notation of \cite[\S 17]{Fu}, this equals $\pi_*(f^*p\cdot [\pi])$, where $[\pi]$ is the orientation class.  Let $\chi_i=e_i^*\tilde{p}$, where $e_i:\eK\to D$ is the $i$-th evaluation map and $\tilde{p}\in N^1(D)$ is the class of a point.

The following relation between these classes leads to the analogue of the Caporaso-Harris formula.

\begin{lemma}
\label{relations}
If $r>3d$, then $$h=\sum_{i=1}^n \varrho_i\chi_i.$$
\end{lemma}

\begin{proof}  First we derive the key consequence of the hypothesis $r>3d$ and in doing so fix some notation.  Let $F:\fC\to\pe$ be a twisted stable map over $\bC$ which lies in $\eK_d(\vrho)$, let $C$ be the coarse curve of $\fC$, and let $g:C\to\bP^2$ be the induced morphism.  Let $x_1,\ldots,x_n\in\fC$ be the marked points and let $y_1,\ldots,y_m$ be the twisted nodes of $\fC$ which lie on at least one component of $\fC$ which maps with positive degree.  By renumbering the markings if necessary, assume that $x_1,\ldots,x_k$ lie on components of $\fC$ which map with positive degree and that $x_{k+1},\ldots,x_n$ lie on components which map with degree $0$.

First note that no twisted node $y_i$ can lie on two components which map with positive degree.  Indeed, since the contact types of $y_i$ on the two components must add to $r$ (since the node is balanced), it would follow that the pullback of $E$ to the normalization of $C$ has degree at least $r$, contradicting the condition $r>3d$.  So each $y_i$ lies on a unique component mapping with positive degree.  Let $\sigma_i$ be the contact type of $y_i$ on this component.  For the same reason, the contact types $\sigma_i$ for all $i$ and $\varrho_j$ for $1\le j\le k$ are equal to the intersection number between $E$ and the component of $C$ at the given point.

Suppose that $A\subseteq\{k+1,\ldots,n\}$ is the set of markings lying in a fixed connected component of $g^{-1}(E)$ and suppose that $B\subseteq\{1,\ldots,m\}$ is the set of nodes lying in the same component.  Note that every irreducible component of $C$ contained in $g^{-1}(E)$ maps with degree $0$.  We claim that 
\begin{equation}
\label{eq:multiplicities}
\sum_{i\in B}\sigma_i=\sum_{j\in A}\varrho_j.
\end{equation}
That they are congruent modulo $r$ follows from Lemma \ref{contact_deg0} together with the balanced condition at the nodes.  From this congruence, equality is deduced from the fact that both sides are between $0$ and $r-1$.

To show that $h=\sum\varrho_i\chi_i$, it suffices to show that for any one-dimensional integral closed substack $V\subseteq\eK_d(\vrho)$, $$\int_V h=\sum\varrho_i\int_V\chi_i.$$  It also suffices to replace $V$ with its normalization.  Let $\eC_V$ be the restriction of the universal curve to $V$.  Then $$\int_V h=\int_{\eC_V} f_V^* p = \int_{\bP^2} p\cap (f_V)_*[\eC_V],$$ which is the degree of $f_V$.  There exists a finite, flat base change $W\to V$ such that every irreducible component of $\eC_W$ has generically irreducible fibers over $W$.  Replace $V$ with $W$.

Now we use the notation introduced at the beginning of the proof, where we take $F:\fC\to\pe$ to be a map corresponding to a general point of $V$.  For each node $y_i$, there are two irreducible components of $\eC_V$ which contain $y_i$, and the intersection between these components is a section $t_i:V\to\eC_V$.  Let $s_i:V\to\eC_V$ be the section corresponding to the $i$-th marking.  The degree of $f_V$ can be computed by $(f_V)_*(f_V)^*[E]=(\deg(f_V))[E]$.  This shows that $$\deg(f_V)=\sum_{i=1}^k\varrho_i\deg(f_V\circ s_i)+\sum_{i=1}^m\sigma_i\deg(f_V\circ t_i).$$  Here the morphisms on the right hand side are viewed as going from $V$ to $E$.  The first summation is $\sum_{i=1}^k\varrho_i\int_V\chi_i$.  It follows from equation \ref{eq:multiplicities} that the second summation is $\sum_{i=k+1}^n\varrho_i\int_V\chi_i$.
\end{proof}

Let $\alpha=(\alpha_1,\ldots)$ and $\beta=(\beta_1,\ldots)$ be sequences with $I\alpha+I\beta=3d$.  Let $\vrho$ be an $|\alpha|+|\beta|$-tuple of integers such that $\alpha_k$ of the numbers $\varrho_1,\ldots,\varrho_{|\alpha|}$ equal $k$ and $\beta_k$ of the remaining numbers are $k$ for each positive integer $k$.  Let $\eK_d^*(\vrho)\subseteq \eK_d(\vrho)$ be the union of the irreducible components where the general source curve is smooth and maps birationally onto its image.  By Theorem \ref{enumerative} together with integration over fibers, we have
\begin{equation}
\label{Nd_int}
N_d(\alpha,\beta)=\frac{1}{\beta!}\int_{\eK_d^*(\vrho)} h^{|\beta|-1}\cdot \prod_{i=1}^{|\alpha|}\chi_i.
\end{equation}

This equation together with Lemma \ref{relations} gives us the following Theorem.

\begin{thm}
\label{CH_recursion}
If $|\beta|>1$, then $$N_d(\alpha,\beta)=\sum_{k:\beta_k>0} kN_d(\alpha+e_k,\beta-e_k).$$
\end{thm}

\begin{proof} If $j$ is chosen so that $\varrho_j=k$ and $j>|\alpha|$, then equation \ref{Nd_int} implies
$$N_d(\alpha+e_k,\beta-e_k)= \frac{\beta_k}{\beta!}\int_{\eK^*_d(\vrho)} h^{|\beta|-2}\cdot\left(\prod_{i=1}^{|\alpha|}\chi_i\right)\cdot\chi_j.$$
This makes the right hand side of the above equation equal to
$$\frac{1}{\beta!}\int_{\eK_d^*(\vrho)} \left(h^{|\beta|-2} \prod_{i=1}^{|\alpha|}\chi_i\right) \sum_{j=|\alpha|+1}^{|\alpha|+|\beta|} \varrho_j\chi_j.$$
To prove that this equals the left hand side, it suffices to show that in $N^1(\eK^*_d(\vrho))$, we have $h=\sum_{i=1}^n \varrho_i\chi_i$ and $\chi_i^2=0$.  The former equality follows from Proposition \ref{relations} and the latter equality is obvious, because if one chooses two distinct points of $E$, then their preimages under the $i$-th evaluation map are disjoint.
\end{proof}

This formula reduces all the numbers of the form $N_d(\alpha,\beta)$ to numbers of the form $N_d(\alpha+\beta-e_k,e_k)$.  We now go a step further, reducing everything to numbers $M_d(\alpha+\beta)$, where $M_d$ is defined as follows.  Let $\gamma$ be a sequence with $I\gamma=3d$ and let $\varphi$ be the following morphism.
$$\xymatrix@R=0pc{
\times_i E^{\gamma_i}\ar[r]^{\varphi} & \Pic^{3d}(E) \\
(x_{ij})_{1\le j\le \gamma_i}\ar@{|->}[r] & \sO_E(\sum ix_{ij})}$$
Let $W=\varphi^{-1}([\sO_E(d)])$, where $[\sO_E(d)]$ denotes the point of $\Pic^{3d}(E)$ corresponding to the restriction of $\sO_{\bP^2}(d)$.  With appropriate choices of base points, $\varphi$ is a homomorphism of abelian varieties and $W$ is the kernel.  Hence the connected components $W_1,\ldots,W_m$ of $W$ are irreducible.  We define $$M_d(\gamma)=\frac{1}{m}\sum_{i=1}^m \deg(ev^{-1}(W_i)\to W_i),$$ where $ev:(K^{\gamma}_d)^*\to W$ is the product of evaluation maps.

We remark that if $n$ is the greatest common divisor of the set of $i$ for which $\gamma_i\not= 0$, then $m=n^2$.  Since we do not need this fact, we omit the proof.  However, the following fact is important.

\begin{lemma}
\label{degree}
Let $$q_k:\times_i E^{\gamma_i}\to \times_i E^{(\gamma-e_k)_i}$$ be the projection which forgets a factor in the $k$-th position.  Then for any connected component $W_i$ of $W$, the degree of $q_k|_{W_i}$ is $k^2/m$.
\end{lemma}

\begin{proof}
Consider the following commutative diagram.
$$\xymatrix{W \ar[r]^{q_k|_W} \ar[d] & E^{|\gamma|-1} \ar[d]^{\mathrm{id}\times [\sO_E(d)]} \\
E^{|\gamma|} \ar[r]_-{q_k\times\varphi} & E^{|\gamma|-1}\times\Pic^{3d}(E)}$$
The morphism $q_k\times\varphi$ has the same degree as the $k$-th power map $\Pic^1(E)\to\Pic^k(E)$, which is $k^2$.  Since both vertical arrows are injective with cokernel $\Pic^{3d}(E)$, the horizontal arrows have the same kernel.  It follows that $\deg(q_k|_W)=k^2$.  Since $q_k|_W$ is a surjective homomorphism and $E^{|\gamma|-1}$ is connected, it follows that each connected component of $W$ maps onto $E^{|\gamma|-1}$ with the same degree.
\end{proof}

\begin{thm}
\label{nice_formula}
For any sequences $\alpha$ and $\beta$ with $I\alpha+I\beta=3d$ and $|\beta|\ge 1$,
$$N_d (\alpha,\beta) = I^\beta \cdot I\beta \cdot \frac{(|\beta|-1) !}{\beta !} \cdot M_d(\alpha+\beta).$$
\end{thm}

\begin{proof}
We prove this by induction on $|\beta|$, starting with $|\beta|=1$.  We must show that $N_d(\alpha,e_k) = k^2M_d(\alpha+e_k)$.  Using the above notation, with $\gamma=\alpha+e_k$, we have the following diagram.
$$\xymatrix{ & W \ar@{^{(}->}[d] & \\
(K_d^{\alpha+e_k})^* \ar[ur] \ar[r]_{ev} & E^{|\alpha|+1} \ar[r]_{q_k} & E^{|\alpha|}}$$
We have shown that $N_d(\alpha,e_k)=\deg(q_k\circ ev)$.  Since $ev$ factors through $W$, we can write this as $$\sum_{i=1}^m \deg(ev^{-1}(W_i)\to W_i)\cdot\deg(q_k|_{W_i}).$$  So the case $|\beta|=1$ follows from Lemma \ref{degree}.

If $|\beta|>1$ we can apply Theorem \ref{CH_recursion}, which by induction gives us
\begin{align*}
N_d(\alpha,\beta) &= \sum_{k:\beta_k>0} k N_d(\alpha+e_k,\beta-e_k) \\
&= \sum_{k:\beta_k>0} k\frac{I^{\beta}}{k}(I\beta-k)\frac{(|\beta|-2)!}{\beta!}\beta_k M_d(\alpha+\beta) \\
&= \frac{I^{\beta}(|\beta|-2)!}{\beta!}\sum_k\beta_k(I\beta-k) M_d(\alpha+\beta)\\
&= \frac{I^{\beta}(|\beta|-2)!}{\beta!}\left(I\beta|\beta|-I\beta\right) M_d(\alpha+\beta) \\
&= I^{\beta}\cdot I\beta\cdot\frac{(|\beta|-1)!}{\beta!}\cdot M_d(\alpha+\beta)
\end{align*}
\end{proof}

{\it Remark.} Since Lemma \ref{relations} holds on the entire stack $\eK_d^{\gamma}$, the above theorem holds for the corresponding Gromov-Witten invariants $I_d(\alpha,\beta),$ provided that one removes the factor of $\beta!$ from the denominator.  Another way to see this is to apply the formulas for multiple cover contributions in Theorem \ref{enumerative}.

\section{How to compute the numbers}

We summarize by showing how one can compute the numbers $N_d(\alpha,\beta)$.  Using Theorem \ref{enumerative}, it is sufficient to compute the corresponding Gromov-Witten invariants $I_d(\alpha,\beta)$.  If there is a $k>1$ such that $\beta_k>0$, then one can apply the WDVV recursion (\ref{WDVV_equation}) to reduce it to invariants which either have lower degree or which have the same degree but have more contacts with $E$ ($|\alpha|+|\beta|$ increases).  Inside of this recursion, one must use Theorem \ref{large_contact_type} to substitute for the invariants involving large contact types.  If there is no such $k$, then one can use Theorem \ref{nice_formula} to express it in terms of an invariant where $|\alpha|$ is smaller (for the case of $|\alpha|+|\beta|=2$, see the remark at the end of section 6).  Eventually, this reduces everything to the numbers $N_d(0,3de_1)$, which can be computed via Kontsevich's recursion.  
This has been implemented in a maple program which can be downloaded from \texttt{http://www.math.lsa.umich.edu/\~{}cdcadman}.

If one could get a sufficient handle on the moduli stacks involved, it might be possible to compute the numbers $N_d(0,e_{3d})$.  An additional obstacle is the fact that some of the components of the moduli stack which contribute have greater than the expected dimension.  Therefore, it would not be easy to compute the contribution to $I_d(0,e_{3d})$ coming from multiple covers and reducible curves.

\bibliographystyle{amsalpha}
\bibliography{cubic}

\end{document}